\documentclass[12pt]{article}

\usepackage{amssymb}
\usepackage{amsmath}
\usepackage{graphicx}
\usepackage[english]{babel}

\newtheorem{theorem}{Theorem}
\begin{document}
\begin{center}
{\MakeUppercase{Exact analytical  solution for one nonlinear variational problem of the cavitation  theory}}

\vspace{0.5em}
{\MakeUppercase{D.V.\,Maklakov and I.R.\,Kayumov}}

\vspace{0.5em}

 \parbox{0.8\linewidth}{In this paper we investigate the limiting
values of the lift and drag coefficients of profiles
in the Helmholtz--Kirchhoff  (infinite
cavity) flow. The coefficients are based on the wetted arc length of
profile surfaces. Namely, for a given value of the lift coefficient
we find minimum and maximum values of the drag
coefficient. Thereby, we  determine
maximum and minimum values of the lift-to-drag ratios.}
\end{center}

\vspace{1.em}


 In the theory of aero and hydrofoils
  there are known two classical models for studying flows past  a profile.
  For the first model the flow is  continuous (fig.\,1a)) and
  for the second one the flow is  separated with formation of an infinite cavity (fig.\,1b)). If we assume that the flow is steady, irrotational and incompressible, then
 for the first model the drag force $D=0$ (d'Alembert's paradox) and the lift
force $L$ is defined by the well-known Kutta--Joukowskii theorem:
\begin{equation}\label{KG}
L=-\rho v_0\Gamma, \quad \Gamma=\int_0^l(\mathbf{v}\cdot
\boldsymbol\tau)d s.
\end{equation}
Here $\rho$ is the density of the fluid, $v_0$ is the velocity at infinity,
  $\Gamma$ is the circulation around the profile,
  $l$ is the perimeter of the profile surface,
   $s$ is the arc abscissa of the profile contour,
  reckoned from the trailing edge point $A$, $(\textbf{v}\cdot
\boldsymbol\tau)$  is the dot product of the velocity vector $\mathbf{v}$
at the point on the profile surface
 and the tangential unit vector $\boldsymbol\tau$, directed toward increase of  $s$.
 For the continuous model the point  $B$ with  the arc abscissa  $s=l$ coincides with the
 point  $A$ for which  $s=0$. If $l_1$ is the arc abscissa of the stagnation point  $O$
and  $v=|\mathbf{v}|$, then
\begin{equation}\label{signs}
(\mathbf{v}\cdot \boldsymbol\tau)=-v(s)\ \
\text{for}
\ \ 0\le s \le l_1,\ \
 (\mathbf{v}\cdot \boldsymbol\tau)=v(s)
\ \
\text{for}
\ \  l_1\le s\le l.
\end{equation}

\begin{figure}[h!]
\centering \includegraphics[width=0.8\linewidth,clip=]{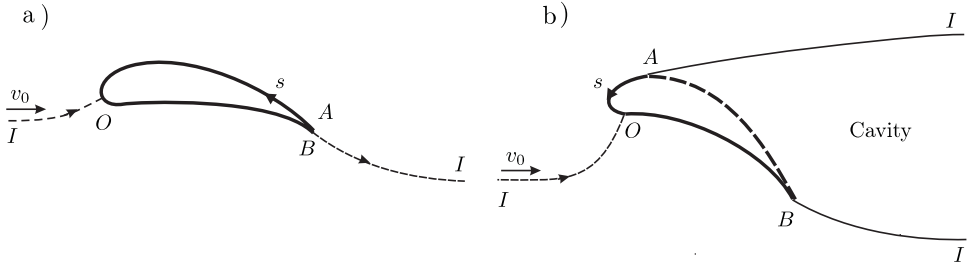}
\caption{a) --- continuous flow over an aerofoil, b)--- Helmholtz--Kirchhoff flow with an infinite cavity past a profile. }
\label{fig1}
\end{figure}

As one  can see from  \eqref{KG}, to compute the lift force for the continuous
model  one needs
only to know the velocity distribution $v(s)$ along the profile surface. Moreover, if $v(s)$
is known, the contour of the profile can be restored by means of solving the so-called
inverse boundary-value problem  of aerodynamics \cite{eip}. The Kutta--Joukowskii theorem
played an outstanding role in the theory of aerofoils and was used many times
for aerodynamic shape optimization (see, for example, \cite{eip,ekm}).

Consider now the second classical model --- the Helmholtz--Kirchhoff  (infinite
cavity) flow (fig.\,1b)).
  According to this model the flow detaches from the profile surface at
the points  $A$ and $B$, and  an infinite  cavity
with a constant pressure, equal to the incident pressure, forms  behind the profile.
The velocity on the free streamlines $AI$ and  $BI$ is constant and equals the incident
velocity $v_0$. As previously
 the stagnation point is denoted  by  $O$   and  the arc abscissa $s$
is reckoned from the point $A$.

For the  Helmholtz--Kirchhoff flow formulas analogous to  \eqref{KG} have
been recently obtained  in the works \cite{makl1}, \cite{makl2}:
\begin{equation}\label{LDmakl}
L=\rho v_0 \int_0^l (\mathbf{v}\cdot\ \boldsymbol\tau)\,
\log\frac{v_0}{v}\, d s, \quad D=\frac{\rho v_0}{4 \pi}
\left(\int_0^l \frac{v}{\sqrt{\varphi}}\, \log\frac{v_0}{v} ds \right)^2,
\end{equation}
where  $l$ is the length of the wetted arc $AOB$ of the profile,
$\varphi=\varphi(s)$ is the distribution of potential
along  $AOB$:
$$
\varphi=\int_s^{l_1} v(s) d s
\ \
\text{for}
\ \ 0\le s \le l_1,\ \
\varphi=\int_{l_1}^s v(s) d s
\ \
\text{for}
\ \ l_1\le s\le l,
$$
$l_1$ is the arc abscissa of the critical point  $O$.

Nowadays the Helmholtz--Kirchhoff model is treated as a limiting case of  cavity
flows (\cite{gur}, p.\,23), when the pressure in the cavity tends to  the incident pressure,
and the sizes of the cavity become infinitely large.
In the theory of cavity flows
the first Brillouin condition \cite{gur} plays an important role: the pressure
in the cavity is minimal. Then the velocity on the free streamlines
$AI$ and   $BI$ is maximal, and therefore
\begin{equation}\label{Brill}
v(s)\le v_0, \quad 0\le s \le l.
\end{equation}
This implies that in formulas  \eqref{LDmakl} the factor
$\log\frac{v_0}{v}\ge 0$.

Due to the simplicity of formulas \eqref{LDmakl} one can formulate
different optimization problems in which one needs
 to determine a velocity distribution
that satisfies the Brillouin condition \eqref{Brill} and  has some optimal property.
One of such problems has been solved  in the works \cite{makl1,makl2}.
Namely, it has been found the velocity distribution that under
the  Brillouin condition \eqref{Brill} provides a global
maximum of the lift force. It has been established that for
the profile of maximum lift the length $l_1=0$ (the points $A$ and  $O$ coincide) and
the optimal velocity distribution $v(s)=e^{-1} v_0=const$, where $e$ is the base of natural logarithms.
It follows from \eqref{LDmakl} that $L_{\max}=\rho v_0^2 l
e^{-1}$  and  this is the global maximum of the lift force.
But such a formulation does not take
into account at all the cavitation drag, which is defined by the second equation
in \eqref{LDmakl}. If $v(s)=e^{-1} v_0$, then
 according to \eqref{LDmakl}  the drag force $D= \rho v_0^2 l/(\pi
e)$, and the lift-to-drag ratio of the profile of maximum lift is
 $\varkappa=L/D=\pi$.
To obtain profiles with a greater lift-to-drag ratio it seems to be natural
to introduce the drag $D$ in the  optimization process.

Let us introduce the lift and drag coefficients $C_L$ and
$C_D$:
$$
C_L=\frac{2 L}{\rho v_0^2 l},\quad C_D=\frac{2 D}{\rho v_0^2 l},
$$
based on the wetted arc length  $l$.

At the end of the paper \cite{makl1} as a variant of a
further perspective direction of investigations it has been formulated the problem
of finding velocity distributions that provide  maximum of the
lift-to-drag ratio $\varkappa=L/D=C_L/C_D$
under the given lift force $L$ (given $C_L$).
 In this paper we present an exact analytical solution
to this problem. Besides, for sake of completeness,  we find  velocity distributions
that provide minimum of the lift-to-drag ratio $\varkappa$.   Thereby, for  fixed values
of the lift coefficient $C_L$ we determine exact upper and lower bounds of this important for applications  hydrodynamic characteristic.

Let  $l_2=l-l_1$ be the length of the arc  $OB$. We introduce two
dimensionless functions  $u_1(\sigma)$ and $u_2(\sigma)$, $0\le \sigma\le 1$,
such that
\begin{equation}\label{u12}
\frac{v}{v_0}=
\begin{cases}
 u_1\big(\frac{l_1-s}{l_1}\big)&\text{on $OA$};\\
 u_2\big(\frac{s-l_1}{l_2}\big)&\text{on $OB$}.
 \end{cases}
\end{equation}
Since the velocity $v\ge 0$, the functions  $u_1(\sigma)$ and $u_2(\sigma)$
are nonnegative. Under the Brillouin condition \eqref{Brill} they satisfy
the inequalities
\begin{equation}\label{Brillu}
u_1(\sigma)\le 1,\quad u_2(\sigma)\le 1.
\end{equation}
By means of  \eqref{LDmakl}  we
 express the lift and drag coefficients  in terms of $u_1(\sigma)$ and
 $u_2(\sigma)$:
  \begin{equation}\label{maklLD}
C_L=2 \{(1-\varepsilon) I[u_2]-\varepsilon I[u_1]\},
\quad C_D=\frac{1}{2\pi}\left\{\sqrt{1-\varepsilon}
J[u_2]+\sqrt{\varepsilon} J[u_1]\right\}^2,
\end{equation}
where  $\varepsilon=l_1/l$, $I[u]$ and  $J[u]$ are nonlinear functionals
of $u(\sigma)$, $0\le\sigma\le 1$:
\begin{equation}\label{maklF}
I[u]=-\int_0^1 u(\sigma)\log u(\sigma)\,d\sigma, \quad
J[u]=-\int_0^1 \frac{u(\sigma) \log u(\sigma)\,d\sigma}
{\sqrt{\int\limits_0^\sigma u(\sigma_1)d \sigma_1}}.
\end{equation}
As one can see from   \eqref{maklF}, under the Brillouin condition
\eqref{Brillu} the values of the functionals~$I[u]$ and $J[u]$ at
$u=u_1(\sigma)$ and $u=u_2(\sigma)$ are nonnegative.

Let us rewrite $I[u]$ and  $J[u]$ in terms of classical
functionals of  calculus of variations. To do so we transform
the function $u(\sigma)$ to $\lambda(\sigma)$:
$$
\lambda(\sigma)=\sqrt{2\int_0^\sigma u(\sigma_1)d\sigma_1}.
$$
Then
\begin{equation}\label{IJ}
I[\lambda]=-\int_0^1 \lambda \lambda' \log ( \lambda \lambda')d
\sigma, \quad J[\lambda]=-\sqrt{2} \int_0^1 \lambda' \log ( \lambda
\lambda') d \sigma.
\end{equation}

It is clear that $\lambda(\sigma)\ge 0$. Besides,
$u(\sigma)=\lambda(\sigma)\lambda'(\sigma)$, whence
$\lambda'(\sigma)\ge 0$.

The  solution of  the basic problem of finding absolute extrema of the lift-to-drag
ratio $\varkappa$ under the given value of $C_L$ is based
on solving the following auxiliary problem.

\textbf{Auxiliary problem.} \textit{Find the function} $\lambda(\sigma)$,
$\sigma\in [0,1]$:
\begin{equation}\label{maincond}
 \lambda(0)=0,\quad \lambda'(\sigma)\ge 0,
\end{equation}
\textit{which delivers a global minimum (maximum) to the functional $J[\lambda]$
under the constraint  $I[\lambda]=q$ ($q$ is given), and the complementary
condition}
\begin{equation}\label{mainaux}
 \lambda(\sigma)\lambda'(\sigma)\le 1.
\end{equation}

Without the complementary condition  \eqref{mainaux} the  auxiliary problem
is a constrained problem of calculus of variations  with a free right endpoint.
Let us find a solution  by the Lagrange multiplier rule without regard for
the nonstandard condition \eqref{mainaux}. To do so we construct
the augmented cost  functional
\begin{equation}\label{extf}
P[\lambda]=-\int_0^1\lambda'(\lambda+k)\log(\lambda\lambda') d
\sigma= -\int_0^1 E(\lambda,\lambda') d \sigma,
\end{equation}
where  $k$ is a real constant.  We write the Euler equation \cite{lavr}
$$
E_\lambda-\frac{d}{d\sigma}E_{\lambda'}=0,
$$
which for the functional  $P[\lambda]$ takes the form
$$
\frac{\lambda''(\lambda+k)}{\lambda'}+\lambda'=0.
$$
Integrating this equation yields
\begin{equation}\label{difur}
\lambda'(\sigma)[\lambda(\sigma)+k]=c,
\end{equation}
where $c$ is a constant. Because  the right endpoint of the desired function
$\lambda(\sigma)$ is free, there holds  the relation  $[E_{\lambda'}]_{\sigma=1}=0$
(see \cite{lavr}), which can be  reduced to
\begin{equation}\label{trance}
\lambda(1)\lambda'(1)=1/e.
\end{equation}
Equation \eqref{difur} subject to the condition $\lambda(0)=0$
 can be easily integrated and has two solutions:
\begin{align}\label{extr}
\lambda(\sigma)=-k+\sqrt{2 c\sigma+k^2}
\ \
\text{for}
\ \ k>0,\\
\label{extr1}
\lambda(\sigma)=-k-\sqrt{2 c\sigma+k^2}
\ \
\text{for}
\ \ k<0.
\end{align}
We demonstrated that the functions of the form  \eqref{extr}, \eqref{extr1}
define  the maximum of the functional $J[\lambda]$, and as a result,
do not define the  minimum of the lift-to-drag ratio  $\varkappa$, but
its maximum, which we  determine in the paper only for sake of completeness.
 Thus, application of the classical approach  does not lead to finding extrema that are of most
practical interest.
To solve  the auxiliary problems we use the technique developed earlier
in papers  \cite{makopt}--\cite{makl_defl} for investigation of
extremal problems of the jet and cavity theory.
This technique is based on Jensen's inequality (\cite{hardy}, theorem~204).

We denote by  $J_{\min}(q)$ and $J_{\max}(q)$, correspondingly,
the global  minimum   and maximum of the functional  $J[\lambda]$ for a given value
of $I[\lambda]=q$. A full solution to the auxiliary problem is given by
\begin{theorem}\label{th1}
$1)$ The  function $J_{\min}(q)$ is defined by the parametric equations
\begin{equation}\label{par1}
\begin{split}
q&=q(b)=\tfrac{1}{2}\big[b^2-k^2+k(b-a)-k^2 \log \tfrac{b}{a}\big],\\
J_{\min}&=\sqrt{2}\big(k \log\tfrac{b}{a}+b+k\big),
\end{split}
\end{equation}
where $b\in (\sqrt{2/e},\sqrt{2})$,
\begin{equation}\label{param1}
k=K(b)=-\frac{(e-1)b(b^2 e -2)}{2+(e-2) b^2 e},\quad
a=\frac{ b^2 e -2}{(e-1) b}.
\end{equation}
The global minimum is achieved by the function
$$
\lambda(\sigma)=
\begin{cases}
\sqrt{2\sigma}&\text{for $0\le \sigma <\gamma$}, \\
-k+\sqrt{2 c(\sigma-\gamma)+(a+k)^2}& \text{for $\gamma\le \sigma \le 1$},
\end{cases}
\quad c=\frac{2-b^2}{2+(e-2) b^2 e}.
$$

$2)$ The function   $J_{\max}(q)$ is defined by the parametric equations
\begin{equation}\label{par2}
\begin{split}
q&=q(b)=\tfrac{1}{2}\big[b^2+k\, b-k^2 \log \tfrac{b+k}{k}\big],\\
J_{\max}&=\sqrt{2}\big(k \log\tfrac{b+k}{k}+b\big),
\end{split}
\quad
\text{where}
\ \ k=K_1(b)=-\frac{b(b^2e-2)}{2(b^2 e-1)},
\end{equation}
             $b\in (0,\sqrt{2/e}]$.

The global maximum is achieved
by the functions
\begin{alignat*}{2}
\lambda(\sigma)&=-k-\sqrt{2 c\sigma+k^2} &\ \ &\text{for}\ \ q\in (0,q_*),\qquad b\in (0,\sqrt{1/e}),\\
\lambda(\sigma)&=\sigma/\sqrt{e} &\ \ & \text{for}\ \ q=q_*,\\
\lambda(\sigma)&=-k+\sqrt{2 c\sigma+k^2} &\ \ & \text{for}\ \ q\in
(q_*,q_{\max}],\ \ b\in (\sqrt{1/e},\sqrt{2/e}],
\end{alignat*}
where  $c=\frac{b^2}{2(b^2 e-1)}$.
\end{theorem}

A full solution to the problem of finding the absolute maximum of the lift-to-drag ratio
$\varkappa$ is given  by
\begin{theorem}\label{th2}
The maximum value of the lift coefficient $C_L=2/e$.  At a given value of the lift coefficient  $C_L\in(0,2/e]$
the global minimum of the drag coefficient is
 $C_{D\min}=\frac{1}{2\pi}J_{min}^2(C_L/2)$,
and the global maximum of the lift-to-drag ratio
 $\varkappa_{\max}=2 \pi C_L/
J_{min}^2(C_L/2)$.
\end{theorem}

The problem
of finding   the  minimum of the
lift-to-drag ratio (the maximum of the drag)  at a given $C_L$   turns out to be
more complex. We established that
\begin{equation}\label{CDmax}
C_{D\max}=\max_{\varepsilon, q_1,q_2}\frac{1}{2
\pi}\left\{\sqrt{1-\varepsilon} J_{\max}(q_2)+\sqrt{\varepsilon}
J_{\max}(q_1)\right\}^2,
\end{equation}
where  $\varepsilon$, $q_1$ and  $q_2$ satisfy the constraints
\begin{equation}\label{restrq}
\varepsilon\in[0,1),\quad 0\le q_1, q_2\le q_{\max}, \quad
(1-\varepsilon) q_2-\varepsilon q_1=C_L/2\in(0,q_{\max}],\quad
q_{\max}=1/e.
\end{equation}
The solution to the problem    \eqref{CDmax}, \eqref{restrq} was found
numerically by means of the standard function  {\bf Maximize} of the package
 Mathematica~8.0.

\begin{figure}[h!]
\centering \includegraphics{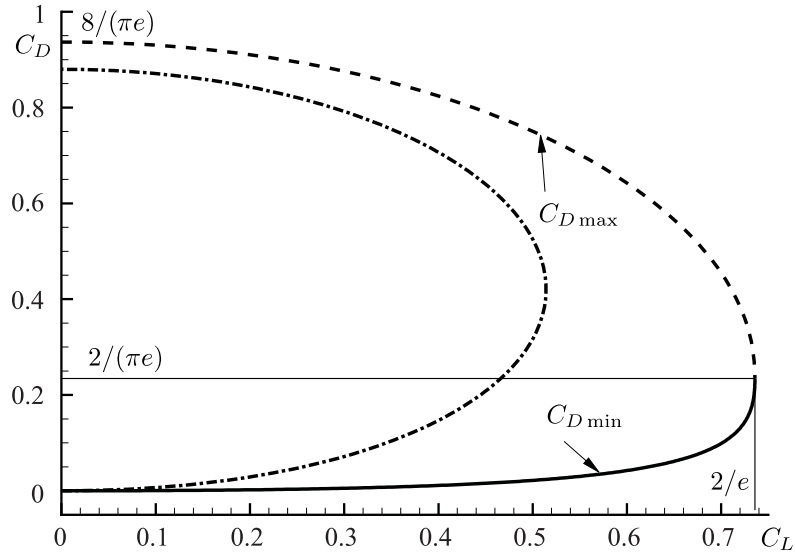}
\caption{Dependencies $C_{D\min}$ and  $C_{D\max}$ on  $C_L$
(the dash-and-dot line is the dependence  $C_D$ on
$C_L$ for a flat plate).}
\label{fig4}
\end{figure}

In fig.~\ref{fig4} we demonstrate the dependencies of the minimal drag coefficient
$C_{D\min}$ and maximal drag coefficient
  $C_{D\max}$ on the lift coefficient $C_L$, and in table~\ref{tb1}
we show  the maximal and minimal lift-to-drag ratios
 $\varkappa_{\max}$,
$\varkappa_{\min}$ for different  $C_L$.
In fig.~\ref{fig4} the dash-and-dot line demonstrates
the dependence  $C_D$ on  $C_L$ for the flat plate.
The coefficients  $C_D$ and  $C_L$ for the flat
plate are defined by Rayleigh's well-known formulas  (\cite{gur}, p.\,83):
$$
C_D=\frac{2\pi \sin^2 \alpha}{4+\pi\sin \alpha}, \quad C_L=\frac{\pi \sin 2\alpha}{4+\pi\sin \alpha},
\quad \varkappa=\cot \alpha,
$$
where $\alpha$ is the angle of attack.

\begin{center}
\begin{table}[h!]
\caption{The values of  $\varkappa_{\max}$ and  $\varkappa_{\min}$ for different
 $C_L$}\label{tb1}
\begin{tabular}{cccccccccc}  
 $C_L$ & 0 & 0.1 & 0.2 & 0.3 & 0.4 & 0.5 & 0.6
  & 0.7 & $\frac{2}{e}$ \\
 $\varkappa_{\max}$ & $\infty$ & 224.88 & 99.1015 &
  57.0649 & 35.9197 & 23.0608 & 14.1997 & 7.0821 &
  $\pi $ \\
 $\varkappa_{\min}$ & 0 & 0.107495 & 0.219695 &
  0.342541 & 0.48536 & 0.666406 & 0.933793 &
  1.53824 & $\pi$
\end{tabular}
\end{table}
\end{center}

As one can see from fig.\ref{fig4} the dash-and-dot line lies entirely
between the curves  $C_{D\min}(C_L)$ and
$C_{D\max}(C_L)$.
It is worthy of  note that for any profile in the Helmholtz--Kirchhoff flow
the point  $(C_L,C_D)$ always lies between the cures
 $C_{D\min}(C_L)$ and  $C_{D\max}(C_L)$.

\end{document}